\documentclass[11pt]{article}
\usepackage[final]{epsfig}
\usepackage{graphics}
\usepackage{amsmath}
\usepackage{amsfonts}
\usepackage{latexsym}
\usepackage{amssymb}
\usepackage{graphicx}
\usepackage{url}
\usepackage{epstopdf}
\usepackage{color}
\usepackage{marginnote}
\usepackage{a4wide}
\usepackage{subcaption}

\newtheorem{lemma}{Lemma}[section]
\newtheorem{proposition}[lemma]{Proposition}
\newtheorem{remark}[lemma]{Remark}
\newtheorem{convention}[lemma]{Convention}
\newtheorem{example}[lemma]{Example}
\newtheorem{theorem}{Theorem}

\newcommand{\g}{{\gamma}}
\newcommand{\G}{{\Gamma}}
\newcommand{\eps}{{\varepsilon}}
\newcommand{\proofend}{$\Box$\bigskip}
\newcommand{\C}{{\mathbb C}}

\newcommand{\R}{{\mathbb R}}
\newcommand{\Z}{{\mathbb Z}}
\newcommand{\RP}{{\mathbb {RP}}}

\newcommand{\A}{\mathcal{A}}
\newcommand{\B}{\mathcal{B}_a}
\newcommand{\Crit}{\mathrm{Crit}}
\def\proof{\paragraph{Proof.}}

\renewcommand{\P}{\mathcal{P}}

\title{Symplectically convex and symplectically star-shaped curves\\ -- \\ a variational problem}

\author{Peter Albers\footnote{
Mathematisches Institut,
Universit\"at Heidelberg,
69120 Heidelberg,
Germany;
peter.albers@uni-heidelberg.de}
 \and 
 Serge Tabachnikov\footnote{
Department of Mathematics,
Pennsylvania State University,
University Park, PA 16802,
USA;
tabachni@math.psu.edu}
} 

\date{\today}

\begin{document}
\maketitle
\begin{abstract}
In this article we propose a generalization of the 2-dimensional notions of convexity resp.~being star-shaped to symplectic vector spaces. We call such curves symplectically convex resp.~symplectically star-shaped. After presenting some basic results we study a family of variational problems for symplectically convex and symplectically star-shaped curves which is motivated by the affine isoperimetric inequality. These variational problems can be reduced back to two dimensions. For a range of the family parameter extremal points of the variational problem are rigid: they are multiply traversed conics. For all family parameters we determine when non-trivial first and second order deformations of conics exist. In the last section we present some conjectures and questions and two galleries created with the help of a Mathematica applet by Gil Bor.
\end{abstract}


\section{Introduction} \label{sec:intro}

In the paper we make a step toward expanding some notions and results of  equiaffine differential geometry of the plane to symplectic spaces. 

Let $\g$ be a smooth, closed, strictly convex, positively oriented plane curve. One can give the curve a parameterization $\g(t)$ such that $[\g'(t),\g''(t)]=1$ for all $t$, where the bracket denoted the determinant made by two vectors. This is called an equiaffine parameterization and, accordingly, one defines the equiaffine length of the curve.

The affine isoperimetric inequality between the equiaffine length $L$ and the enclosed area $A$ assert
$$
L^3\leq 8\pi^2 A,
$$
with equality if and only if $\g$ is an ellipse, see, e.g., \cite{Lut,ST}.  Note that the inequality goes in the ``wrong'' direction, compared to the usual isoperimetric inequality!

Assume that the origin is inside the curve, then $\g$ is also star-shaped, in addition to being convex. For any parameterization $\g(s)$, one has a well-defined (i.e., independent of the parameterization) differential 1-form and a cubic form
$$
[\g(s),\g'(s)]\ ds,\quad  [\g'(s),\g''(s)]\ ds^3,
$$
such that
$$
A=\frac{1}{2} \int [\g(s),\g'(s)]\ ds,\quad L=\int \sqrt[3]{[\g'(s),\g''(s)]}\ ds.
$$
Thus the affine isoperimetric inequality relates the integrals of these two 1-forms along a convex closed curve.

Let $\g(t)$ be a smooth closed curve in the standard symplectic vector space $(\R^{2n},\omega)$. Call $\g$ \emph{symplectically star-shaped} if $\omega(\g(t),\g'(t))>0$ for all $t$, and \emph{symplectically convex} if $\omega(\g'(t),\g''(t))>0$ for all $t$. 

\begin{remark}
{\rm
We point out that an alternative definition of symplectically star-shaped resp.~convex is to require $\neq$ instead of $>$ above. That this is actually more general is explained in Section \ref{subsec:sconv} where we construct examples of curves with all possible sign combinations. 
}
\end{remark}

Similarly to the plane case, we may define two differential 1-forms along a symplectically star-shaped and symplectically convex curve $\g$ by
$$
 \omega(\g,\g')\ dt,\quad \sqrt[3]{\omega(\g',\g'')}\ dt.
$$
Both forms are well defined and have no zeroes. Inspired by the affine isoperimetric inequality, we are interested in the relative extrema of $\int_\gamma \sqrt[3]{\omega(\g',\g'')}\ dt$ constrained by $\int_\gamma \omega(\g,\g')\ dt$. In fact, we consider a more general variational problem: describe the curves $\g$ that are the relative extrema of $\int_\gamma \omega(\g',\g'')^a\ dt$ constrained by $\int_\gamma \omega(\g,\g')\ dt$ where $a$ is a real exponent. We point out that the case $a=\tfrac13$ corresponds to the affine isoperimetric inequality.

Before we come to our main results we point out that the corresponding metric problem, extremizing the $L^2$-norm of the (metric) curvature on a class of plane curves, is a widely studied topic going back to Bernoulli and Euler and goes under the name of \emph{elastica}, see \cite{Sin}. In addition, we mention the recent article \cite{Wan} in which the $L^p$-norms of the curvature are studied and \cite{OPW} where the corresponding gradient flows are developed. The latter is a natural next step also for the affine context from this article. See \cite{ST} for an affine analog of the Euclidean curve shortening flow.

\paragraph{Main results}

We prove that such extremal curves of this variational problem lie in symplectic affine 2-planes, and therefore the problem reduces to a 2-dimensional one (Proposition \ref{prop:planar}).

We then fix the constraint by giving the curve the centroaffine parameterization, i.e.,~we assume $[\g(t),\g'(t)]=1$. Then Hill's equation $\g''(t)=-p(t)\g(t)$ holds, and we may consider  the functional $\B(\g):=\int [\g'(t),\g''(t)]^a\ dt$.  

For $a\in[\tfrac12,1]$, one has rigidity: the extremal curves are multiply traversed conics (Propositions \ref{prop:diffeq} and \ref{lm:int}). The same rigidity result holds, although for a different reason, in the case of the affine isoperimetric inequality, $a=\tfrac13$, (Theorem \ref{thm:1/3}).

In Theorem \ref{thm:infdef} we describe non-trivial infinitesimal deformations of multiply traversed conics in the class of extremal curves: if $a=\tfrac13$, then the $n$-fold ellipse is infinitesimally rigid; otherwise,
a non-trivial infinitesimal deformation of the $n$-fold ellipse exists if and only if
$$
a=\frac{k^2-2n^2}{k^2-4n^2}
$$
for some positive integer $k\neq n$.

Theorem \ref{thm:signdef} concerns second order deformations of conics: for $a<0$, the circle $\g_0$ is a local minimum of $\B$; for $a\in (0,\frac{1}{3})$, it is a local maximum; for $a>\frac{7}{5}$, it is a local minimum; and in other cases the Hessian is not sign-definite. The Hessian is degenerate (with 1-dimensional kernel) if and only if $$a=\frac{k^2-2}{k^2-4}$$ for some positive integer $k$.

In Section \ref{sect:gal} we present  examples of extremal curves and formulate some conjectures about them.

This introduction would not be complete if we failed to mention another reason for our interest in centroaffine differential geometry, namely its close relation with the Korteweg-de Vries equation, discovered by U. Pinkall \cite{Pin} and studied by a number of authors since then. When the exponent $a$ equals $2$, the extremal curves are periodic solutions to Lam\'e's equation thoroughly studied in this context in a recent paper \cite{BBT}.  
\bigskip 

{\bf Acknowledgements}: We are very grateful to Gil Bor for writing a Mathematica program that made it possible to experiment with extremal curves. 
PA was supported by Deutsche Forschungsgemeinschaft (DFG, German Research Foundation) through Germany’s Excellence Strategy EXC-2181/1 - 390900948 (the Heidelberg STRUCTURES Excellence Cluster), the Transregional Colloborative Research Center CRC/TRR 191 (281071066) and the Research Training Group RTG 2229 (281869850).  ST was supported by NSF grant DMS-2005444.

\section{Examples of symplectically convex and symplectically star-shaped curves} \label{subsec:sconv}

In this section we construct curves $\g$ with all possible sign combinations of the quantities $\omega(\g,\g')\neq0$ and $\omega(\g',\g'')\neq0$. In particular, we assume that all curves are immersed. We start with a remark concerning  symplectically star-shaped curves. The sphere $S^{2n-1}\subset \R^{2n}$ carries a  contact structure defined by the symplectic orthogonal complement to the position vector. A symplectically star-shaped curve projects to a  transverse curve in $S^{2n-1}$. A similar remark applies to the contact $\RP^{2n-1}$, the projectivization of $\R^{2n}$. If  $\omega(\g,\g')>0$ then $\g$ is positively transverse and $<$ corresponds to negatively transverse. A somewhat similar interpretation for the condition $\omega(\g',\g'')\neq0$ is derived in Lemma \ref{lm:infl} in case of $\R^4$.

We consider the unit sphere $S^3\subset \C^2=\R^4$ with its standard contact structure. The standard contact form at a point $q\in S^3$ is $\omega(q,\cdot)|_{T_qS^3}$. Let $\g(t)$ be a smooth closed Legendrian curve in $S^3$, i.e., $\omega(\g,\g')=0$. Then $J\g'$ is a vector normal to $\g$ inside the contact plane. Here $J$ is the complex structure on $\C^2$. The following lemma is well known, see, e.g., \cite{Ben}.

\begin{lemma} \label{lm:push}
Pushing a closed Legendrian curve $\g$ slightly in the direction of $J\g'$, resp. $-J\g'$, inside $S^3$ yields a negative, resp. positive, transverse curve. Thus, a Legendrian curve in $S^3\subset \C^2$ can be deformed into a symplectically star-shaped curve.
\end{lemma}

\proof
Let us assume that $\g$ is parametrized by arc-length and set $\Gamma=\g+\eps J\g'$. Then
$$
\omega(\Gamma,\Gamma')=\eps \big[\omega(\g,J\g'')+\omega(J\g',\g')\big] + O(\eps^2).
$$
If we denote by $\cdot$ the inner product then $\omega(J\g',\g')=-\g'\cdot\g' = -1$. From $\g\cdot\g=1$ we conclude $\g\cdot \g'=0$, i.e., $\omega(\g,J\g')=0$. Differentiating this equality gives
$\omega(\g,J\g'') + \omega(\g',J\g')=0$ and hence $\omega(\g,J\g'')=\omega(J\g',\g')$. It follows that for sufficiently small $\eps >0$ one has
$\omega(\Gamma,\Gamma') <0$. That is, $\Gamma$ is a negative transverse curve. The case of $\Gamma=\g-\eps J\g'$ is similar.
\proofend

Let us say that a Legendrian curve $\g$ in $S^3$ has an \textit{inflection point} at $p\in\g$ if it is second-order tangent to its tangent great circle at $p$. A generic curve in $S^3$ does not have inflection points but, as we shall see, a generic Legendrian curve has finitely many of them.

The tangent Gauss map sends a Legendrian curve in $S^3$ to the space of oriented Legendrian great circles, that is, to the oriented Lagrangian Grassmannian $\Lambda_2^+$. The image of the Gauss map is a smooth curve which has vanishing differential precisely  at the points corresponding to inflection points of the Legendrian curve. 

\begin{lemma} \label{lm:infl}
A necessary and sufficient condition for a Legendrian curve $\g$ having an inflection point at $\g(t)$ is $\omega(\g'(t),\g''(t)) = 0$.
\end{lemma} 

\proof
The condition $\omega(\g',\g'') = 0$ is invariant under reparameterization of $\g$, so we may assume that $\g$ is parameterized by arc-length. 

First, we claim that the orthogonal projection of  $\g''$ to $S^3$ is $\g+\g''$. Indeed,  $\g\cdot \g'=0$ implies that $\g'\cdot \g' + \g\cdot \g''=0$, hence $\g\cdot \g''=-1$. Now $\g\cdot (\g+\g'')=1-1=0$, as claimed. Therefore, an inflection point is characterized by the vanishing of $\g+\g''$.

Next, we claim that the tangential acceleration vector $\g+\g''$ lies in the contact structure. Indeed, $\omega(\g,\g')=0$, hence, after differentiating, $\omega(\g,\g'')=0$. Therefore, $\omega(\g,\g+\g'')=0$, as needed. 

In addition, $\g\cdot\g=\g'\cdot\g'=1$ implies that $\gamma'$ and $\g+\g''$ are orthogonal to each other. Finally, since $\g$ is Legendrian, i.e.,~$\omega(\g,\g')=0$, the tangent vector $\g'$ lies also in the contact structure. To summarize, we have two orthogonal vectors, $\g'$ and $\g+\g''$, in the contact structure. In particular, $\omega(\g',\g+\g'')=\pm\|\g'\| \|\g+\g'\|=\pm\|\g+\g'\|$ since $\|\g'\|=1$.

The Lemma now follows  from $\omega(\g',\g'')=\omega(\g',\g+\g'')=0$ if and only if $\g+\g''=0$, i.e.,~if and only if $\g$ has an inflection point.
\proofend

It follows that if we find a closed Legendrian curve $\g$ with $\omega(\g',\g'') \neq 0$, then, according to Lemma \ref{lm:push} its small push in the normal direction in the contact plane will yield a curve with  $\omega(\g,\g') \neq 0$ and $\omega(\g',\g'') \neq 0$. Let us show how to construct such a Legendrian curve $\g$. See \cite{Ar,OT} for related matters.

Consider the Reeb field of the standard contact form on $S^3$ and let $pr: S^3\to S^2$ be the respective Hopf fibration. The projection $p$ takes the Legendrian great circles in $S^3$ to the great circles in $S^2$. Closed Legendrian curves are projected to closed  immersed curves in $S^2$ that bound a region with signed area that is a multiple of $2\pi$ and, conversely, such a spherical curve lifts (non-uniquely) to a closed Legendrian curve in $S^3$. 

Since $pr$ is a Riemannian submersion it follows that the inflections of a Legendrian curve $\g$ correspond to the inflections of the spherical curve $pr(\g)$. Accordingly, every closed spherical curve with everywhere positive geodesic curvature and bounding area $2\pi k$ for some $k\in \Z$ lifts to a Legendrian curve which is free from inflections. Starting with a closed spherical curve with everywhere positive geodesic curvature, the area condition is easily arranged by adding appropriately sized kinks to the curve. 

Finally, we note that the  map $(z_1,z_2)\mapsto (\bar z_1,\bar z_2)$ of $\C^2$ preserves the contact structure and the property of $\g$ being Legendrian but changes the sign of $\omega(\g'(t),\g''(t))$ to the opposite. Therefore we can have all four combinations of signs of the quantities $\omega(\Gamma,\Gamma')$ and $\omega(\Gamma',\Gamma'')$ where $\Gamma$ is a small push-off as in Lemma \ref{lm:push}.

\begin{remark}
{\rm
Another interpretation of the above construction is by looking at a different Hopf fibration $\widetilde{pr}: S^3\to S^2$, whose fibers are Legendrian. This Hopf fibration takes Legendrian great circles to circles (of some, possibly zero radius) on $S^2$. The projection of a smooth Legendrian curve $\g$ is a wave front, possibly with cusps, and the inflections of $\g$ correspond to the vertices of the spherical curve $\widetilde{pr}(\g)$. If $\widetilde{pr}(\g)$ is smooth and convex, it must have at least four vertices (the 4-vertex theorem), without convexity at least two vertices, but if $\widetilde{pr}(\g)$ has cusps it may be vertex-free. 
}
\end{remark}


\section{Towards the solution of the variational problem} \label{subsec:var}

In this section we describe the setting of our variational problem and prove that extremizers are contained in symplectic affine 2-planes. This allows us to reduce the problem to the 2-dimensional case.

The problem is to find the closed symplectically star-shaped and symplectically convex curves that are extrema of $\B(\gamma):=\int\omega(\g',\g'')^{a}\ dt$ , $a\in \R\setminus\{0\}$, subject to the constraint given by $\A(\gamma):=\int\omega(\g,\g')\ dt$. More precisely, we consider the space $ \P_T$ of $T$-periodic symplectically star-shaped and symplectically convex curves in $\R^{2n}$ and consider $\A,\B:\P_T\to\R$ and ask for extrema of $\B$ subject to $\A=c_0$, that is, we want to describe the set 
$$
\Crit(\B|_{\{\A=c_0\}})\subset \P_T.
$$ 
Let us note that both functionals, $\A$ and $\B$, are translation invariant, that is, do not depend on the choice of the origin (as long as the curve remains star-shaped). We start with a few simple observations.

\begin{lemma}
For $c_0>0$ and $c_1>0$ there is a natural bijection (by rescaling) from $\P_T$ to itself inducing a bijections
$$
\{\A=c_0\}\cong\{\A=c_1\}
$$
and
$$
\Crit(\B|_{\{\A=c_0\}})\cong\Crit(\B|_{\{\A=c_1\}}).
$$ 
\end{lemma}

\proof
Let $\gamma_0\in\P_T$ and consider
$$
\gamma_1(t):=\left(\tfrac{c_1}{c_0}\right)^{\frac12}\gamma_0(t)\in\P_T.
$$
If $\gamma_0\in\{\A=c_0\}$ then
$$
\A(\gamma_1)=\tfrac{c_1}{c_0}\A(\gamma_0)=c_1,
$$
i.e., $\gamma_1\in\{\A=c_1\}$. From
\begin{equation}\label{eq:scaling_of_B}
\B\left(\left(\tfrac{c_1}{c_0}\right)^{\frac12}\gamma\right)=\left(\tfrac{c_1}{c_0}\right)^a\B(\gamma)
\end{equation}
it follows that the bijection $\gamma\mapsto \left(\tfrac{c_1}{c_0}\right)^{\frac12}\gamma$ just rescales $\B$ by a fixed factor, i.e., induces the claimed bijection $\Crit(\B|_{\{\A=c_0\}})\cong\Crit(\B|_{\{\A=c_1\}})$.
\proofend

\begin{remark}
{\rm Equation \eqref{eq:scaling_of_B} implies that, if $\B$ has critical points at all, they appear in $\R_{>0}$-family and the critical value is necessarily $0$.
}
\end{remark}

\begin{lemma}
The bijection $\P_T\to\P_1$ given by $\Gamma(t):=\gamma(t/T)$ preserves $\A$ and rescales $\B$ by $\tfrac{1}{T^{3a-1}}$.
\end{lemma}

\proof
We consider the bijection $\P_T\to\P_1$ given by $\Gamma(t):=\gamma(t/T)$. Then
\begin{equation}\nonumber
\begin{aligned}
\A(\Gamma)&=\int_0^1\omega(\Gamma(t),\Gamma'(t))\ dt=\int_0^1\omega\big(\gamma(\tfrac{t}{T}),\tfrac1T\gamma'(\tfrac{t}{T})\big)\ dt=\int_0^T\omega\big(\gamma(s),\gamma'(s)\big)\ ds=\A(\gamma),
\end{aligned}
\end{equation}
as it has to be since $\A(\Gamma)=\A(\gamma)$ is twice the enclosed area. Similarly, we compute
\begin{equation}\nonumber
\begin{aligned}
\B(\Gamma)&=\int_0^1\omega(\Gamma'(t),\Gamma''(t))^a\ dt\\
&=\int_0^1\tfrac{1}{T^{3a}}\omega\big(\gamma'(\tfrac{t}{T}),\gamma''(\tfrac{t}{T})\big)^a\ dt=\tfrac{1}{T^{3a-1}}\int_0^T\omega\big(\gamma'(s),\gamma''(s)\big)^a\ ds\\
&=\tfrac{1}{T^{3a-1}}\B(\gamma),
\end{aligned}
\end{equation}
as claimed. \proofend

Now, we begin to study the relative extrema of $\B$ constrained by $\A$. We call the extrema \emph{critical curves}. The previous two lemmata imply that we may consider curves of fixed period and with fixed constraint. 

\begin{proposition} \label{prop:planar}
A critical curve is contained inside a symplectic affine 2-plane of $(\R^{2n},\omega)$.
\end{proposition}

\proof 
We first derive the equation for a critical curve using a Lagrange multiplier $\lambda\in\R$. For that let $v(t)$ be a vector field along $\g$ and consider an infinitesimal deformation $\g_\eps=\g+\eps v$ of $\g$. Then $\g$ is critical if for some $\lambda$ for all $v$ 
$$
\frac{d}{d\eps}\bigg|_{\epsilon=0} \left[\lambda \int \omega(\g_\eps,\g_\eps')\ dt + \int \omega(\g_\eps',\g_\eps'')^{a}\ dt   \right] =0.
$$
The derivative of the first integral in direction of $v$ computes to 
$$
\int (\omega(v,\g')+\omega(\g,v'))\ dt = 2 \int \omega(v,\g')\ dt,
$$
where we used integration by parts and skew-symmetry of the symplectic form in the equality sign. 

Similarly, setting $F:=\omega(\g',\g'')^{a-1}$, the derivative of the second integral can be expressed as 
\begin{equation*}
\begin{split}
a \int F \big(\omega(v',\g'')+\omega(\g',v'')\big) \ dt 
= a \int \big(F' \omega(v',\g')+2F\omega(v',\g'')\big)\ dt \\
= a \int \big(F'' \omega(\g',v)+3F' \omega(\g'',v)+2F \omega(\g''',v) \big)\ dt.  
\end{split}
\end{equation*}
Since $v$ is arbitrary and $\omega$ is non-degenerate, the criticality condition implies that the vectors $\g',\g'',\g'''$ are linearly dependent, and since $F\neq 0$, one has
$
\g'''=f\g'+g\g''
$
for some periodic functions $f(t),g(t)$.

It follows that the bivector $\g'\wedge\g''$ satisfies the differential equation $(\g'\wedge\g'')'=g(\g'\wedge\g'')$, that is, remains proportional to itself. Hence the curve $\g'$ is planar. It follows by integration that the curve $\g$ lies in an affine 2-plane, parallel to the plane spanned by $\g'$ and $\g''$.
\proofend

Proposition \ref{prop:planar} allows us to reduce the problem to two dimensions as follows. A critical curve $\g$ is contained in an affine 2-plane. Let $V$ be the two dimensional linear space parallel to this affine 2-plane, and $W$ be its symplectic orthogonal. Then we may write $\g=\g_1\oplus \g_2 \in V \oplus W$.  Since $\g'\in V$ we conclude that $\g_2$ is constant. We claim that $\g_1$ is also a critical curve. Indeed, since $W$ is symplectically orthogonal to $V$ we compute
 $$
 0\neq \omega(\g,\g')=\omega(\g_1\oplus \g_2,\g_1') = \omega(\g_1,\g_1'),
 $$
i.e., $\g$ is symplectically star-shaped and $\A(\g)=\A(\g_1)$. Trivially, $\g_1$ is also symplectically convex with $\B(\g)=\B(\g_1)$. Now, that we have reduced the problem to two dimensions,  we denote the symplectic form $\omega$ in the plane by brackets $[\cdot,\cdot]$. 
\begin{convention}
To fix the constraint $\A$, from now on we always parametrize a star-shaped curve so that $[\g,\g']=1$, that is, $\g$ is in its centroaffine parameterization. 
\end{convention}
In centroaffine parameterization we have $\g''=-p\g$. The function $p$ is called the \emph{centroaffine curvature} of the curve $\g$. For some computations below we record that $p=[\g',\g'']>0$ and $\g'''=-p'\g-p\g'$.

Thus, we can reformulate our goal as describing all extremal curves of the functional $\B(\g):=\int p^a(t) dt$ on the space of periodic curves in centroaffine parameterization with $p>0$. Here $a\in \R\setminus\{0\}$ is fixed. To be quite precise, we point out that we, of course, consider the space of all periodic curves in centroaffine parameterization of a \emph{fixed period}, mostly, $2\pi$ or a multiple of $2\pi$. As explained in the introduction, the case $a=\tfrac13$ is of special interest since this corresponds to the affine isoperimetric inequality. 

\begin{lemma}\label{lm:invariance}
The properties of a curve being symplectically star-shaped and convex are invariant under point-wise multiplication by an element in $SL(2,\R)$ and sufficiently small translation. For any value of $a$, the functional $\B$ is invariant under this $SL(2,\R)$ action. In addition, if $a=\tfrac13$, the functional $\mathcal{B}_\frac13$ is invariant under translations. This is the only value of $a$ with this additional symmetry.
\end{lemma}

\proof
Point-wise multiplication of a curve with a matrix from $SL(2,\R)$ does not change its centroaffine parameterization nor its centroaffine curvature. 

A sufficiently small translation of a curve $\g$ keeps the properties of being symplectically star-shaped and convex but in general not the centroaffine parameterization of $\g$. After reparametrization, the shifted curve $\bar\g=\g+c$ is in centroaffine parameterization, i.e.~$\bar\g(\tau)=\g(t(\tau))+c$ and $[\bar\g,\frac{d\bar\g}{d\tau}]=1$. Its centroaffine curvature $\bar p$ satisfies 
$$
\bar p (\tau) = \left(\frac{dt}{d\tau}\right)^3  p(\tau),
$$
indeed, 
$$
\bar p=\left[\frac{d\bar\g}{d\tau},\frac{d^2\bar\g}{d\tau^2}\right]=\left[\frac{dt}{d\tau}\frac{d\bar\g}{dt}, \left(\frac{dt}{d\tau}\right)^2\frac{d^2\bar\g}{dt^2} + \frac{d^2t}{d\tau^2}\frac{d\bar\g}{dt} \right]=\left(\frac{dt}{d\tau}\right)^3\left[\frac{d\bar\g}{dt}, \frac{d^2\bar\g}{dt^2} \right]=\left(\frac{dt}{d\tau}\right)^3p
$$
since $\frac{d\bar\g}{dt}=\frac{d\g}{dt}=\g'$. We conclude that precisely for $a=\tfrac13$ the functional $\B$ satisfies 
$$
\B(\bar\g)=\int\bar p^a(\tau)d\tau=\int p^a(t)dt=\B(\g)
$$
for all symplectically star-shaped and convex curves $\g$.
\proofend

We now derive a critical point equation for $\B$ in terms of the function $F=p^{a-1}$.
\begin{proposition} \label{prop:diffeq}
For $a\neq1$ the extremal curves $\g$ of the functional $\int p^a(t) dt$ satisfy  
\begin{equation}\label{eq:difF}
F'''=-2(2+b)F^b F'
\end{equation}
where $F=p^{a-1}$ and $b=\tfrac{1}{a-1}$. If $a=1$ or $a=\tfrac12$, then $p$ is constant and $\g$ is therefore a conic. 
\end{proposition}

\proof
First, we describe the vector fields $v$ along $\g$ that preserve the centroaffine parameterization of $\g$. Write such a field as $v=g\g+f\g'$. Since the deformation by $v$ is assumed to preserve the centroaffine parameterization, we conclude that $[\g,v']+[v,\g']=0$ and hence $2g+f'=0$. Thus, the vector field $v$ has the form $-\frac{f'}{2}\g+f\g'$; the deformations keeping the centroaffine parameterization depend on one periodic function.

As in the proof of Proposition \ref{prop:planar}, linearizing $\int p^a(t)dt$ in $v$ leads to 
$$
\int F\big([\g',v'']+[v',\g'']\big)\ dt =0
$$
for every vector field $v$ as above, where we recall that $F=p^{a-1}$. 
Using integration by parts twice and recalling that $\g'''=-p'\g-p\g'$, we rewrite the integral as
$$
\int \Big(\big(F''-2pF\big)[\g',v]-\big(3pF'+2p'F\big)[\g,v]\Big)\ dt.
$$
Using $[\g,v]=f$, $[\g',v]=\tfrac12f'$ and another integration by parts, this integral becomes
$$
\int \Big(-\tfrac12F'''-2pF'-p'F\Big) f\ dt.
$$ 
The critically condition is that this integral vanishes for all $f$, and we conclude that the integrand is zero.

If $a=1$, then $F\equiv1$, and hence the criticality condition simply becomes $p'=0$. Therefore, $\g$ solves $\g''=\text{const}\cdot\gamma$, i.e., $\g$ is a conic. 

Otherwise, recall that $p=F^b$, hence $p'=b F^{b-1}F'$. Substitute this to the integrand and collect terms to obtain the claimed formula \eqref{eq:difF}.

If $a=\tfrac12$ then $b=-2$, and equation \eqref{eq:difF} reduces to $F'''=0$. Since $F$ is periodic and positive, $F$ is necessarily a constant, and so is $p=F^{-2}$. Thus, again $\g$ is a conic in this case.
\proofend

\begin{proposition} \label{lm:int}
For $a\in(\tfrac12,1)$, equation \eqref{eq:difF} has only constant solutions.
\end{proposition}


\proof
We can integrate equation \eqref{eq:difF} to
\begin{equation}\label{eq:difF_integrated}
F''=-\frac{2(b+2)}{b+1}F^{b+1}+c
\end{equation}
with some constant $c$. Note that $b\neq-1$. Since $F$ is a positive and periodic function, we get at the minimum $m:=\min F$ and the maximum $M:=\max F$ of $F$ the usual inequalities for $F''$, leading to
$$
0\geq -\frac{2(b+2)}{b+1}M^{b+1}+c\qquad\text{and}\qquad0\leq-\frac{2(b+2)}{b+1}m^{b+1}+c.
$$
Thus, we arrive at the two inequalities
\begin{equation}\label{eq:inequalities_F''}
\left\{
\begin{aligned}
&\;\;\frac{b+2}{b+1}m^{b+1}\leq \frac{b+2}{b+1}M^{b+1}\\[.5ex]
&\;\;m\leq M.
\end{aligned}
\right.
\end{equation}
If $b+1>0$ or $b+2>0>b+1$, then the first inequality in (\ref{eq:inequalities_F''}) is consistent with the second. However, the first inequality together with $0>b+2$ implies $m\geq M$ and we conclude $m=M$. That is, $0>b+2$ implies that $F$ is constant. Now, recall that $b=\tfrac{1}{a-1}$. Therefore, $F$ is constant if $\tfrac12<a<1$, as claimed.
%
%
%
%
%
%
\proofend

\begin{lemma}
For $a\not\in[\tfrac12,1]$, equation \eqref{eq:difF} has non-constant positive, periodic solutions.
\end{lemma}

\begin{remark}
{\rm The solutions we construct in the proof of the lemma below are, by construction, often times multiply covered.}
\end{remark}

\proof 
As a preparation, we consider the Hamiltonian formulation of the ODE \eqref{eq:difF_integrated}. Hamilton's equations for the Hamiltonian function ($c\in\R$ is a constant which we will choose appropriately below)
\begin{equation}\nonumber
H(Q,P):=\frac12P^2+\frac{2}{b+1}Q^{b+2}-cQ:\R^2\longrightarrow\R
\end{equation}
are given by
\begin{equation}\nonumber
\left\{
\begin{aligned}
\;\; \dot P &= -\frac{\partial H}{\partial Q} = -\frac{2(b+2)}{b+1}Q^{b+1}+c\\
\dot Q &= \frac{\partial H}{\partial P} = P,
\end{aligned}
\right.
\end{equation}
and are, with $Q=F$ and $P=F'$, equivalent to the ODE \eqref{eq:difF_integrated}. We compute
$$
dH(Q,P)= P\,dQ +\left(\frac{2(b+2)}{b+1}Q^{b+1}-c\right)dQ
$$
and
$$
\mathrm{Hess}H(Q,P)=
\begin{pmatrix}
2(b+2)Q^b & 0\\
0 & 1
\end{pmatrix}.
$$
Therefore, $(Q,P)$ is a critical point of $H$ if and only if
\begin{equation}\nonumber
\left\{
\begin{aligned}
P=0\\
\;\; \frac{2(b+2)}{b+1}Q^{b+1} &=c,
\end{aligned}
\right. 
\end{equation}
and then is a local minimum if  
$$
(b+2)Q^b >0.
$$
Our assumption $a\not\in(\tfrac12,1)$ is equivalent to $b\in(-2,-1)\cup(-1,0)\cup(0,\infty)$, (since $b=\tfrac{1}{a-1}$ and $a\neq0$). 

We are searching for non-constant positive, periodic solutions of equation \eqref{eq:difF}. I.e.~we are looking for non-constant periodic orbits of $H$ with $Q=F>0$.

\underline{In case $b\in(0,\infty)$} we choose $c\in\R$ very large and positive. Then the equation $\frac{2(b+2)}{b+1}Q^{b+1}=c$ determines a local minimum of $F$ at, say, $(Q_0,P_0=0)$ with $Q_0$ large and positive. Therefore, the linearized dynamics given by 
$$
\mathrm{Hess}H(Q_0,0)=
\begin{pmatrix}
2(b+2)Q_0^b & 0\\
0 & 1
\end{pmatrix}.
$$
is a very fast rotation of $\R^2$. That is, the linearized dynamics is periodic with very small periodic. Since $(Q_0,0)$ is a local minimum, any level set $\{H=E\}$ with $E$ slightly larger than the energy  $H(Q_0,0)=\frac{2}{b+1}Q_0^{b+2}-cQ_0$ of $(Q_0,0)$ consists of a small circle (and potentially other connected components). This small circle is then a periodic orbit with period approximately that of the linearized dynamics. Since $b\neq0$, the Hamiltonian $H$ is not quadratic and therefore varying the energy value $E$ also changes the period of the periodic orbit in the level set $\{H=E\}$ near $(Q_0,0)$. I.e.~by iterating and varying the level set we can arrange any period we wish.

\underline{In case $b\in(-1,0)$} we choose $c\in\R$ positive and small. Then the equation $\frac{2(b+2)}{b+1}Q^{b+1}=c$ determines a local minimum of $F$ at $(Q_0,P_0=0)$ with $Q_0$ small and positive. The linearized dynamics is again given by 
$$
\mathrm{Hess}H(Q_0,0)=
\begin{pmatrix}
2(b+2)Q_0^b & 0\\
0 & 1
\end{pmatrix}.
$$
This is still a very fast rotation of $\R^2$ since $b<0$, $b+2>0$ and $Q_0>0$ is small. The argument proceeds as in the previous case.

\underline{In case $b\in(-2,-1)$} we choose $c\in \R^2$ negative and large in absolute value. The equation $\frac{2(b+2)}{b+1}Q^{b+1}=c$ determines a local minimum of $F$ at $(Q_0,P_0=0)$ with $Q_0$ small and positive.
The linearized dynamics is yet again a very fast rotation of $\R^2$ since $b<0$, $b+2>0$ and $Q_0>0$ is small. The argument proceeds as in the first case.
\proofend

%

\begin{remark} \label{rmk:elliptic}
{\rm Equation (\ref{eq:difF}) can be further integrated, which is in terms of the Hamiltonian formulation just expressing the preservation of $H$ along solutions. 
\begin{equation} \label{eq:1stord}
(F')^2=-\frac{4}{b+1}F^{b+2}+2cF+d.
\end{equation}
For $a=2$, resp.~$a=\tfrac32$, we have $b+2=3$, resp.~$b+2=4$, and hence $F$ is the Weierstrass elliptic function. 
Since $p=F^b$, in the first case $p$ is also  the Weierstrass function, and in the second case $p$ is its square.
In the former case, the respective equation $\g''=-p\g$ is called the Lam\'e equation, and it was thoroughly studied, see, e.g., 
\cite{WW}.

}
\end{remark}

\section{Infinitesimal rigidity of multiple conics as critical curves} \label{subsec:infdef}

As a multiple conic we take the unit circle traversed $n$ times. This is a critical curve of the functional $\B(\g)=\int_0^{2\pi n} p^a dt$, and we ask whether it admits a non-trivial infinitesimal deformation in the class of critical curves. 

The functional $\B$ is invariant under the action of $SL(2,\R)$; the respective deformations comprise a 3-dimensional space, and we consider them as trivial. In addition, if $a=\tfrac13$, the functional is  invariant under parallel translations, see Lemma \ref{lm:invariance}. In this case, we add this 2-dimensional space to the deformations that we consider as trivial. The $n$-fold circle is infinitesimally rigid if it does not admit non-trivial 
infinitesimal deformations in the class of critical curves.

Let $\g_0(t)=(\cos t,\sin t)$ be the unit circle traversed $n$ times, and let $\g_1=-\frac{f'}{2}\g_0+f\g_0'$ be a vector field along it that defines its infinitesimal deformation. We assume that the period is $2\pi n$, and accordingly, $f(t)$ is a $2\pi n$-periodic function. 

The trivial deformations are described in the following lemma.

\begin{lemma} \label{lm:triv}
The infinitesimal action of $SL(2,\R)$ corresponds to the functions $f$ that are linear combinations of $1, \cos (2t), \sin (2t)$, and parallel translations to the functions $f$ that are linear combinations of $\cos t, \sin t$. 
\end{lemma}

\proof
The deformed curve is $\g_0+\eps\g_1$.
The case $f=1$ corresponds to the rotation of the unit circle.

Let $f(t)=\sin (2t)$. In this case, computing up to $\eps^2$, the curve $\g_0+\eps \g_1$ is the ellipse $(1+2\eps)x^2+(1-2\eps)y^2=1$,
and likewise for $f=\cos (2t)$.

If $f(t)=2\sin t$ then, again mod $\eps^2$,  the curve $\g_0+\eps \g_1$ is the unit circle $(x+\eps)^2+y^2=1$, 
and likewise for $f=2\cos t$.
\proofend

We are ready to describe the infinitesimal rigidity of the circle.

\begin{theorem} \label{thm:infdef}
If $a=\tfrac13$, then  the $n$-fold unit circle is infinitesimally rigid. Otherwise,
a non-trivial infinitesimal deformation of the $n$-fold unit circle exists if and only if
$$
a=\frac{k^2-2n^2}{k^2-4n^2}
$$
for some positive integer $k\neq n$.
\end{theorem}

\proof
The calculations below are modulo $\eps^2$.

Let $\G=\g_0+\eps\g_1$. We have $\G''=-(p_0+\eps p_1) \G$,  hence 
$$
p_0+\eps p_1=[\G',\G'']=[\g_0'+\eps\g_1',\g_0''+\eps\g_1'']=1+\eps ([\g'_0,\g_1'']+[\g_1',\g_0'']).
$$
Therefore
$$
F=(p_0+\eps p_1)^{a-1}=1+\eps (a-1)([\g'_0,\g_1'']+[\g_1',\g_0'']).
$$

We calculate 
$$
\g_1'=-\left(f+\frac{1}{2}f''\right)\g_0+\frac{1}{2}f'\g_0',\ \g_1''=-\left(\frac{3}{2}f'+\frac{1}{2}f'''\right)\g_0-f\g_0',
$$
hence
$$
q:=[\g'_0,\g_1'']+[\g_1',\g_0'']=2f'+\frac{1}{2}f'''.
$$
Compare with \cite{Pin}, where the Korteweg-de Vries equation is interpreted as a flow on centroaffine curves.

The case of $a=1$ was considered earlier (the only solution of the variational problem is a conic); therefore we assume that $a\neq 1$. 

Since the perturbed curve is critical, equation (\ref{eq:difF_integrated}) holds. Write the constant in this equation as $c_0+\eps c_1$. Then (\ref{eq:difF_integrated}) implies
$$
q''=-2(b+2) q +\frac{c_1}{a-1}.
$$
Since $q''$ and $q$ have zero average, we conclude that $c_1=0$. Therefore 
$q''=-2(b+2) q.$

This equation has periodic solutions only when $b+2>0$, and then $q(t)$ is a linear combination of 
$\cos (\sqrt{2(b+2)}t)$ and $\sin (\sqrt{2(b+2)}t)$. For $q$ to be $2\pi n$-periodic, one has to have $\sqrt{2(b+2)}=\tfrac{k}{n}$, that is, 
\begin{equation} \label{eq:quant}
b=\frac{k^2-4n^2}{2n^2} \ \ {\rm or}\ \ a=\frac{k^2-2n^2}{k^2-4n^2}
\end{equation}
for positive integers $k,n$.  Note that since $b\neq 0$, one has $k\neq 2n$.

Let us show that if condition (\ref{eq:quant}) holds, the desired infinitesimal deformations of a multiple circle exists. 

We find $f(t)$ from the equation
$$
2f'+\frac{1}{2}f''' = A\cos\left(\frac{kt}{n}\right) + B\sin\left(\frac{kt}{n}\right).
$$
It follows that, modulo the kernel of the differential operator $\frac12d^3+2d$, the function $f$ is a linear combination of $\cos\left(\frac{kt}{n}\right)$ and $\sin\left(\frac{kt}{n}\right)$. The kernel of the differential operator contributes trivial deformations, and we end up with a 2-dimensional space of deformations.

It remains to see when these deformations are trivial. Since $k\neq 0$ and $k\neq 2n$, the only ``suspicious" case is $k=n$. In this case, $a=\tfrac13$, and indeed, the first harmonics give  trivial deformations, corresponding to parallel translations.  
\proofend

\begin{remark} 
{\rm We note that, in agreement with Lemma \ref{lm:int}, 
$$
a=\frac{k^2-2n^2}{k^2-4n^2}
$$
does not take values in $[\frac{1}{2},1]$: if $k>2n$, then $a>1$, and if $k<2n$, then $a<\tfrac12$. 
}
\end{remark}

\section{Second order deformations of conics} \label{subsec:2ord}

Here we investigate how the functional $\B(\g) = \int p(t)^a dt$ changes under a second order deformation of the unit circle. We recall that the functional is $SL(2,\R)$-invariant, see Lemma \ref{lm:triv}. 

We assume that the period is $2\pi$ and that the curves have the centroaffine parameterization $[\g,\g']=1$. Let $\G=\g_0+\eps\g_1+\eps^2\g_2$, ignoring the higher order terms in $\eps$. As before, $\g_0(t)=(\cos t,\sin t)$, hence $\g_0''=-\g_0$, and $\g_1=-\frac{f'}{2}\g_0+f\g_0'$, where $f(t)$ is a $2\pi$-periodic function. 

The condition $[\G,\G']=1$ implies
\begin{equation} \label{eq:twocond}
 [\g_0,\g_2']+[\g_1,\g_1']+[\g_2,\g_0']=0.
\end{equation}

We have; $\G''=-P\G$, hence 
$$
P=[\G',\G'']=[\g_0'+\eps\g_1'+\eps^2\g_2',\g_0''+\eps\g_1''+\eps^2\g_2'']=1+\eps q(t) + \eps^2 r(t),
$$
where
$$
q=[\g'_0,\g_1'']+[\g_1',\g_0''],\ r=[\g_0',\g_2'']+[\g_1',\g_1'']+[\g_2',\g_0''].
$$
Then
$$
P^a=1+\eps aq + \eps^2 a\left(r+\frac{a-1}{2} q^2\right).
$$
As we already know, 
$$
q=[\g'_0,\g_1'']+[\g_1',\g_0'']=2f'+\frac{1}{2}f'''.
$$

We need to calculate
$$
\int_0^{2\pi} \left(r+\frac{a-1}{2} q^2\right) dt.
$$
Integrating by parts and using $\g_0'''=-\g_0'$, we get
$$
\int r dt = \int ([\g_0',\g_2'']+[\g_1',\g_1'']+[\g_2',\g_0'']) dt = \int (2[\g_2,\g_0']+[\g_1',\g_1'']) dt.
$$
Integrate  equation (\ref{eq:twocond}):
$$
0 = \int ([\g_0,\g_2']+[\g_1,\g_1']+[\g_2,\g_0']) dt = \int (2[\g_2,\g_0']+[\g_1,\g_1']) dt,
$$
therefore
$$
\int r dt  = \int ([\g_1',\g_1'']-[\g_1,\g_1']) dt.
$$
We calculate
$$
[\g_1,\g_1']=f^2+\frac{1}{2}ff''-\frac{1}{4} f'^2,\ [\g_1',\g_1'']=f^2+\frac{1}{2}ff''+\frac{3}{4}f'^2+\frac{1}{4}f'f''',
$$
hence
$$
\int r dt  = \int \left(f'^2+ \frac{1}{4}f'f'''\right) dt = \int \left(f'^2- \frac{1}{4}f''^2\right) dt.
$$

Next,
$$
\int q^2 dt = \int \left(2f'+\frac{1}{2}f'''\right)^2 dt = \int \left(4f'^2-2f''^2+\frac{1}{4}f'''^2\right) dt.
$$

In the case of most interest, $a=\tfrac13$, and we obtain
\begin{equation} \label{eq:third}
\int \left(r+\frac{a-1}{2} q^2\right) dt = \int \left(r-\frac{1}{3} q^2\right) dt = -\frac{1}{12} \int (4f'^2-5f''^2+f'''^2) dt.
\end{equation}

\begin{lemma} \label{lm:max}
The integral in (\ref{eq:third}) is non-negative, and it equals zero if and only if $f$ is a first harmonic.
\end{lemma}

\proof
Let
$$
f'(t) = \sum_k c_k e^{ikt},\ c_{-k}=\bar c_k
$$
be the Fourier expansion. Then 
$$
\int_0^{2\pi} (4f'^2-5f''^2+f'''^2) dt = 2 \sum_{k>0} (4-5k^2+k^4) |c_k|^2.
$$
We have $4-5k^2+k^4 = (k^2-1)(k^2-4)\ge 0$, and the sum is positive unless the only non-zero term is for $k=1$, that is, $f$ is a first harmonic.
\proofend

Now consider the general case:
$$
\int \left(r+\frac{a-1}{2} q^2\right) dt = \int \left((2a-1)f'^2- \frac{4a-3}{4} f''^2+ \frac{a-1}{8} f'''^2 \right) dt.
$$
In terms of the Fourier coefficients of $f'$, this is
$$
2 \sum_{k>0} \left(2a-1-  \frac{4a-3}{4} k^2 + \frac{a-1}{8} k^4 \right) |c_k|^2.
$$
The expression in the parentheses equals
\begin{equation} \label{eq:Fc}
\frac{1}{8} (k^2-4) \left[(a-1)k^2-2(2a-1)\right].
\end{equation}
We also note that the quadratic term of $P^a$ contains the factor $a$.

\begin{theorem} \label{thm:signdef}
For $a<0$, the circle $\g_0$ is a local minimum of $\B$; for $a\in (0,\frac{1}{3})$, it is a local maximum; for $a>\frac{7}{5}$, it is a local minimum, and in other cases the Hessian is not sign-definite. The Hessian is degenerate (with 1-dimensional kernel) if and only if $a=\frac{k^2-2}{k^2-4}$ for some positive integer $k$.
\end{theorem}


\proof
Let $a=1$. Then the sign of (\ref{eq:Fc}) is that of $-(k^2-4)$, which is positive for $k=1$ and negative for $k\geq 3$.

Let $a>1$. Then the sign of (\ref{eq:Fc}) is positive for sufficiently large $k$. The Hessian is positive-definite if this sign is positive for all $k$. When $k=1$, the first factor in (\ref{eq:Fc}) is negative, and so is the second one: $1-3a$. When $k\ge 3$, the first factor is positive, and the second one is positive if and only if $9(a-1)-2(2a-1) > 0$, that is, $a>\tfrac75$.

Let $0 < a < 1$. Then the sign of (\ref{eq:Fc}) is negative for sufficiently large $k$. The Hessian is negative-definite if this sign is negative for all $k$. When $k=1$, the first factor in (\ref{eq:Fc}) is negative, and the second factor equals $1-3a$. Thus (\ref{eq:Fc}) is negative if and only if $a <\tfrac13$. If this inequality is satisfied then, for $k\ge 3$, the second factor in (\ref{eq:Fc}) is $2-4a-(1-a)k^2 < 0$. 

If $a<0$, then the analysis of the preceding paragraph still holds, but the factor $a$ of the quadratic term in $P^a$ changes the sign to the opposite. 

Finally, the Hessian is degenerate when (\ref{eq:Fc}) is zero for some $k$. Solving this for $a$ yields the last claim of the theorem. 
\proofend

The last statement of the theorem agrees with Theorem \ref{thm:infdef}.
The numbers $a=\frac{k^2-2}{k^2-4}$ form the sequence
$$
\frac{1}{3},\ \frac{7}{5},\ \frac{7}{6},\ \frac{23}{21}, \ldots
$$
that converges to $1$.  Each time that $a$ crosses an element of this ``spectrum", the signature of the Hessian changes by 1.

\section{The case of $a=\tfrac13$} \label{subsec:1/3}

We recall that we attempt to describe extremal curves of the functional $\B=\int p^a(t) dt$ on the space of periodic curves in centroaffine parameterization with $p>0$. The case $a=\tfrac13$ corresponds to the affine isoperimetric inequality. In particular, the functional then is translation invariant, see Lemma \ref{lm:invariance}. In this case, $b=-\tfrac32$, and equation (\ref{eq:difF_integrated}) reads
\begin{equation} \label{eq:2ndord}
F''=2F^{-\tfrac12}+c.
\end{equation}
Since $F''\le 0$ at the maximum and $F>0$, we conclude that $c<0$. 

Next, integrate  equation (\ref{eq:2ndord}) to
\begin{equation} \label{eq:nextint}
(F')^2=8F^{\tfrac12}+2cF+d
\end{equation}
(this is the equation of a level curve of the Hamiltonian, see Section \ref{subsec:var}).

Let $F(t)=G^2(t)$ where $G(t)$ is also a positive periodic function. Then (\ref{eq:nextint}) becomes
\begin{equation}    \label{eq:dintG}
(G G')^2 = \frac{c}{2} G^2 + 2G + d
\end{equation}
(again renaming the constants). 

The right hand side of (\ref{eq:dintG}) is a quadratic polynomial in $G$, and it has at least two roots because $G$ is a periodic function that attains maximum and minimum. Since a quadratic polynomial has at most two roots, $G$ oscillates between its maximum and minimum, and has no other critical values. 

\begin{example} \label{ex:circle}
{\rm 
Let us examine the case of a parallel translated $n$ times traversed circle, which is a critical curve:
$$
\g=(A+r\cos\alpha,B+r\sin\alpha),
$$
where $\alpha(t)$ is a function of the centroaffine parameter $t$. We assume that the range of $t$ is $[0,2\pi]$, that of $\alpha$  is $[0,2\pi n]$, and the radius of the circle is $r$. 

The range of the centroaffine parameter is twice the (algebraic) area bounded by the curve, hence $r=n^{-\tfrac12}$. 

We calculate: 
$$
[\g,\g_\alpha]=r^2+rA\cos\alpha+rB\sin\alpha,
$$
 and since $[\g,\g_t]=1$, we have 
$$
\frac{dt}{d\alpha} = r^2+rA\cos\alpha+rB\sin\alpha.
$$
 This can be integrated:
$$
t = r^2 \alpha +rA\sin\alpha-rB\cos\alpha +C,
$$
or
$$
\alpha+A\sin\alpha+B\cos\alpha=nt+C',
$$
with the constants renamed. Since 
$$
A\sin\alpha+B\cos\alpha= -\sqrt{A^2+B^2}\cos(\alpha+\theta)
$$ 
with 
$$
\sin\theta=\frac{A}{\sqrt{A^2+B^2}},\ \cos\theta=-\frac{B}{\sqrt{A^2+B^2}},
$$
we can change the parameter $\alpha$ to obtain a simplified equation
\begin{equation} \label{eq:alphat}
\alpha(t) - A\cos\alpha(t)=nt+C
\end{equation}
(once again renaming the constants)
Then $d\alpha/dt=n(1+A\sin\alpha)^{-1}.$

Next,
$$
p=[\g_t,\g_{tt}]=[\g_\alpha,\g_{\alpha\alpha}] \left(\frac{d\alpha}{dt}\right)^3 = n^2(1+A\sin\alpha)^{-3}.
$$
Since $G=p^{-\tfrac13}$, we conclude that 
\begin{equation} \label{eq:Gt}
G(t)=  n^{-\tfrac23}(1 + A\sin\alpha(t)).
\end{equation}

}
\end{example}

Let us continue with the general case.
Let $0<m<M$ be the minimum and the maximum of $G$, and let $c=-2k^2$. Write the right hand side of  (\ref{eq:dintG}) as $k^2 (G-m)(M-G)$, then the differential equation becomes
\begin{equation} \label{eq:difG}
G G' = k \sqrt{(G-m)(M-G)},
\end{equation}
where we allow the square root to take values at its positive and negative branches.

Set 
$$
\mu=\frac{M+m}{2}, \ \eps = \frac{M-m}{2}.
$$
Note that since $G>0$, we have $\eps < \mu$.

Since $G$ oscillates between $m$ and $M$, let us make another substitution:
\begin{equation} \label{eq:ans}
G(t)=\mu + \eps \sin \varphi (t),
\end{equation} 
where $\varphi(t)$ is not necessarily a periodic function anymore. Since $G$ is $2\pi$-periodic, $\varphi(t+2\pi)=\varphi(t)+2\pi n$ where $n$ is an integer.

Substitute (\ref{eq:ans}) into (\ref{eq:difG}) to obtain $\varphi'(\mu + \eps \sin \varphi)=k$. This differential equation integrates to
\begin{equation} \label{eq:H}
\mu \varphi(t) - \eps \cos \varphi(t) = kt+C.
\end{equation}
Since $0<\eps<\mu$,  the left hand side is a monotonic function of $\varphi$, therefore this functional equation uniquely determines the function $\varphi(t)$.

 Substituting $\varphi(t+2\pi)$ in (\ref{eq:H}), we find that $\mu (\varphi(t+2\pi)-\varphi(t))=2\pi k$.
Since $\varphi(t+2\pi)=\varphi(t)+2\pi n$ with $n\in \Z$, we have $k=\mu n$. Thus we have
\begin{equation} \label{eq:HH}
\varphi(t) - \eps \cos \varphi(t) = nt+C,
\end{equation}
where, as before, we renamed the constants.
 
\begin{theorem} \label{thm:1/3}
If $a=\tfrac13$ then the relative extrema of the functional $\B$ constrained by $\A$ are multiply traversed conics.
\end{theorem}

\proof
Let $\g(t)$ be a $2\pi$-periodic critical curve. Then equation  (\ref{eq:HH}) holds for some $n$, defining  function $\varphi(t)$. 
Observe that equation (\ref{eq:HH}) is identical to equation  (\ref{eq:alphat}) from Example \ref{ex:circle}. Therefore $\varphi(t)$ coincides with the function $\alpha(t)$. Similarly, equations (\ref{eq:ans}) and (\ref{eq:Gt}) coincide, and so function $G(t)$ coincides with that  from Example \ref{ex:circle}.

It follows that the centroaffine curvature $p(t)$ of the curve $\g(t)$ is the same as that of the parallel translated $n$-fold circle $\g_0(t)$. Consider the contact element $(\g(0),\g'(0))$. Acting on $\g_0(t)$ by an element of $SL(2,\R)$, we can arrange for $(\g_0(0),\g_0'(0))$
to coincide with $(\g(0),\g'(0))$. The action of $SL(2,\R)$ does not change the centroaffine curvature, hence the two curves, $\g$ and $\g_0$, satisfy the same second order differential equation and have the same initial data. Therefore they coincide.
\proofend

\begin{remark} {\rm
The Lambert $W$ function is the inverse function of the complex function $z=we^w$, see \cite{La}. The function $\varphi$ defined by equation (\ref{eq:HH}) is related to the Lambert function in the following way. 

Let us assume that $C=0, n=1$ in (\ref{eq:HH}).
Consider the complex function  given by the equation $\xi=\eta - \eps e^{i\eta}$. If $\eta$ is real then  $\mathfrak{Re}\ \xi = \eta - \eps \cos\eta$, the expression that defines the function $\varphi(t)$. 

Let $z=we^w$, and set $z=-i\eps e^{i\xi},\ w=-i\eps e^{i\eta}.$  Then $\ln z = w+\ln w$, that is,  $\xi=\eta - \eps e^{i\eta}$. Therefore the inverse function $\eta(\xi)$ is conjugated to the Lambert function $w(z)$ by the exponential function. 
}
\end{remark}

\section{Pictures and open problems} \label{sect:gal}

\begin{figure}[hbtp]
\centering
\includegraphics[width=5.5in]{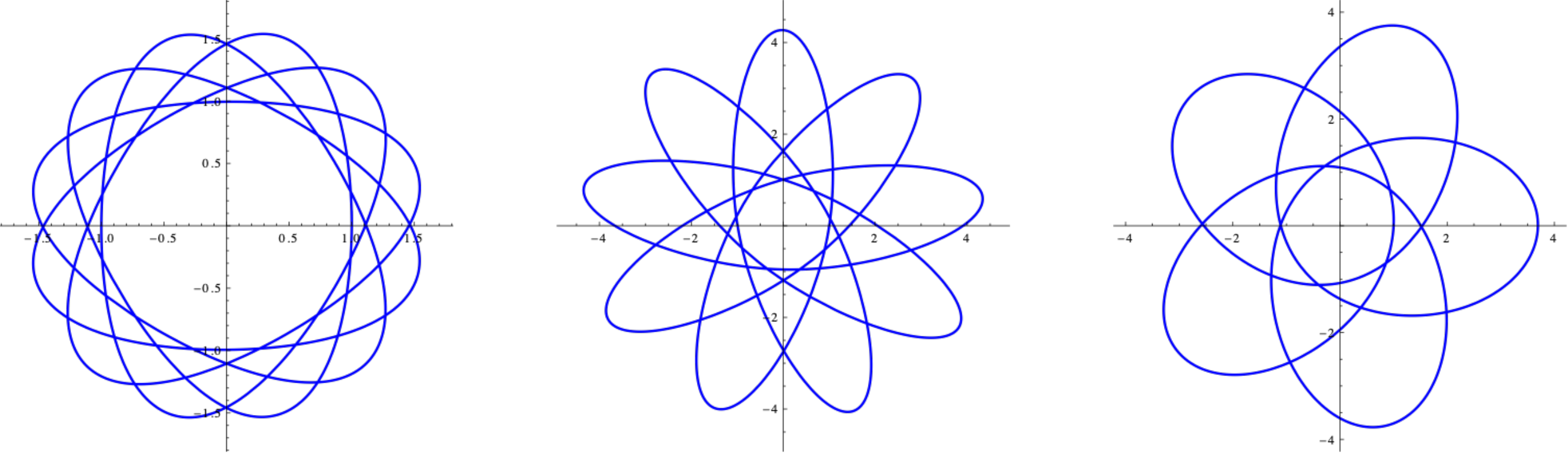}
\caption{Curves having $a=-2, a=-1, a=0.2$ and the rotation numbers $5,5,4$, respectively.}\label{fig:1}
\end{figure}

\begin{figure}[hbtp]
\centering
\includegraphics[width=5.5in]{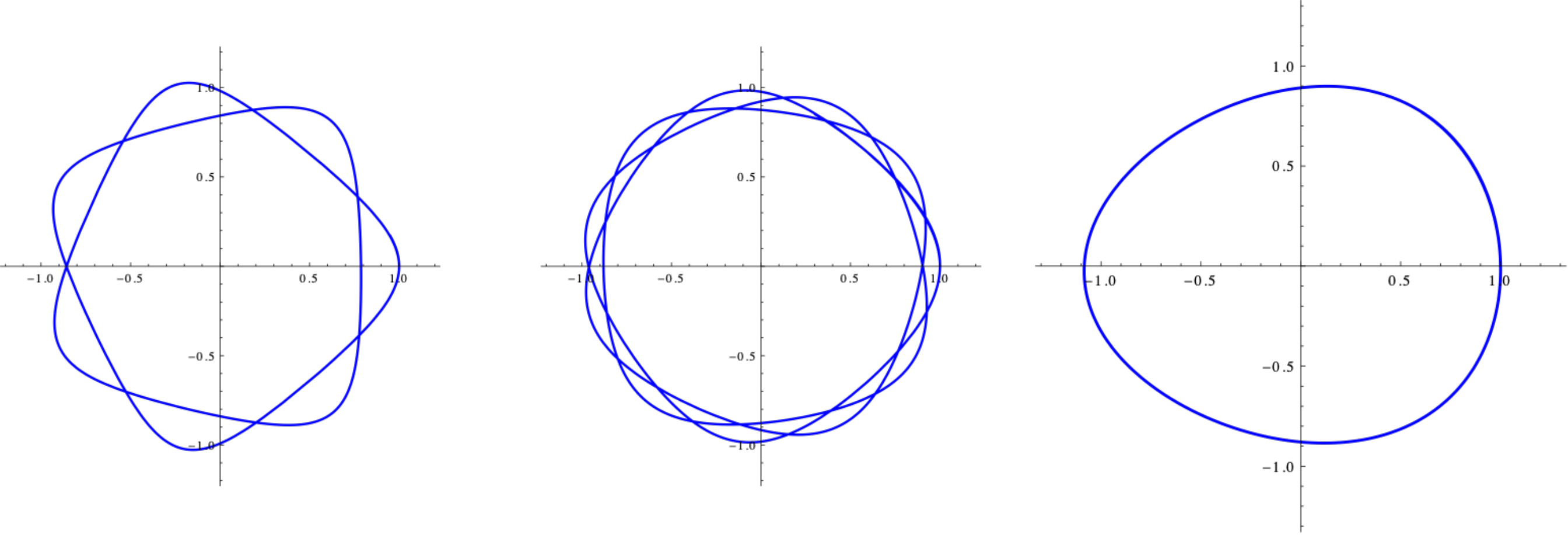}
\caption{These curves have $a=1.2, a=1.2, a=1.4$, respectively. }\label{fig:egg}
\end{figure}

\begin{figure}[hbtp]
\centering
\includegraphics[width=5.5in]{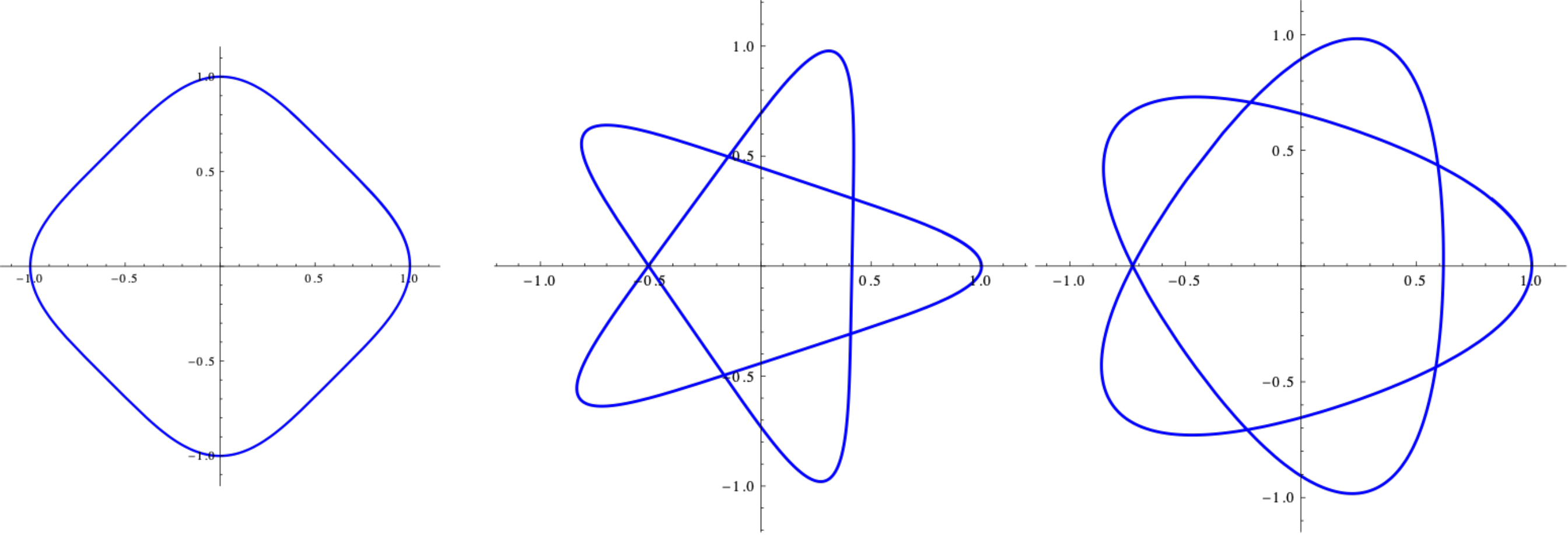}
\caption{These curves have $a=1.5, a=1.5, a=1.75$, respectively. }\label{fig:3}
\end{figure}

\begin{figure}[hbtp]
\centering
\includegraphics[width=5.5in]{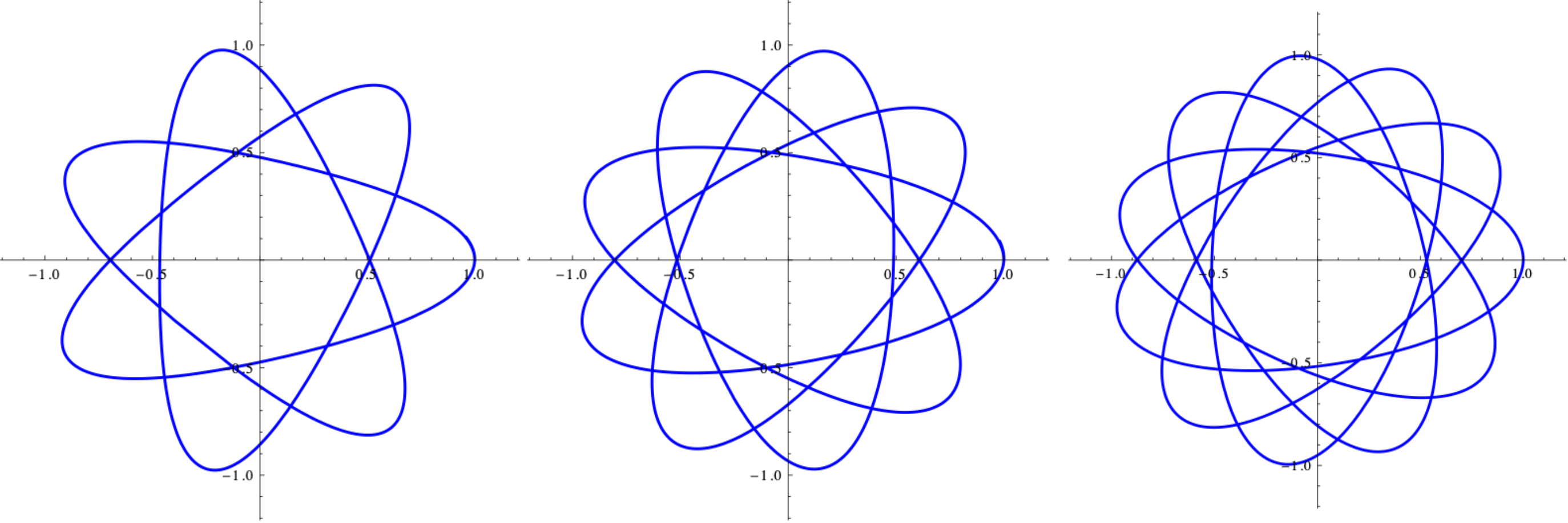}
\caption{These curves have $a=2, a=2.5, a=3$ and the rotation numbers $3,4,5$, respectively.}\label{fig:4}
\end{figure}

First, in Figures \ref{fig:1} -- \ref{fig:4} we present  of extremal curves obtained in numerical experiments using a Mathematica applet created by Gil Bor. 

\clearpage

Let us comment on a common  geometrical feature of these curves. Recall the notion of the osculating circle of a smooth curve in Euclidean geometry: this is a circle that is 2-order tangent to a curve at a given point, that is, it shares the curvature with the curve. Informally speaking, the osculating circle passes through three  infinitesimally close points of the curve.

In centroaffine geometry, the role of osculating circles is played by osculating central conics. The space of central conics is 3-dimensional, and for  every point of a star-shaped curve there exists a central conic that is 2-order tangent to it at this point.  Central conics have constant centroaffine curvature, and the osculating central conic at point $\gamma(t_0)$ has the constant centroaffine curvature $p(t_0)$.

Similarly to the case of osculating circles the osculating central conic goes from one side of the curve to the other side if $p'(t_0) \neq 0$. If $p$ has a non-degenerate maximum or minimum at point $t_0$, that is, $\gamma(t_0)$ is a centroaffine vertex, then the osculating central conic lies on one side of the curve near this point.

In the following lemma we prove that the centroaffine curvature of \emph{extremal curves} takes only two values at their centroaffine vertices, its maximum and minimum.

\begin{lemma} \label{lm:osc}
The function $F$, introduced in Section \ref{subsec:var}, and hence the centroaffine curvature $p$, has only two critical values, its maximum and minimum.
\end{lemma}

\proof
Since the curve is closed, the function $F$ attains its maximum and minimum. It has no other critical values because the right hand side of equation (\ref{eq:1stord}), as a function of $F$, has no more roots than the number of sign changes among its three coefficients, see \cite{PS} (part 5, chapter 1, \S6, No 77).  
\proofend

Looking again at the above figures it is fairly obvious that the case $a=1.4$ in Figure \ref{fig:egg}, the ``egg'' lacks the same symmetry all other curves have. This seems related to Theorem \ref{thm:infdef} asserting that for all $a$ of the form $a=\frac{k^2-2}{k^2-4}$, $k\in\Z_{>2}$, a circle (which is a critical curve for any $a$) admits a non-trivial infinitesimal deformation. Indeed, $1.4=\frac{3^2-2}{3^2-4}$. 

\begin{figure}[hbtp]
\begin{subfigure}{.5\textwidth}
\centering
\end{subfigure}
\begin{subfigure}{.5\textwidth}
\centering
\includegraphics[width=2.5in]{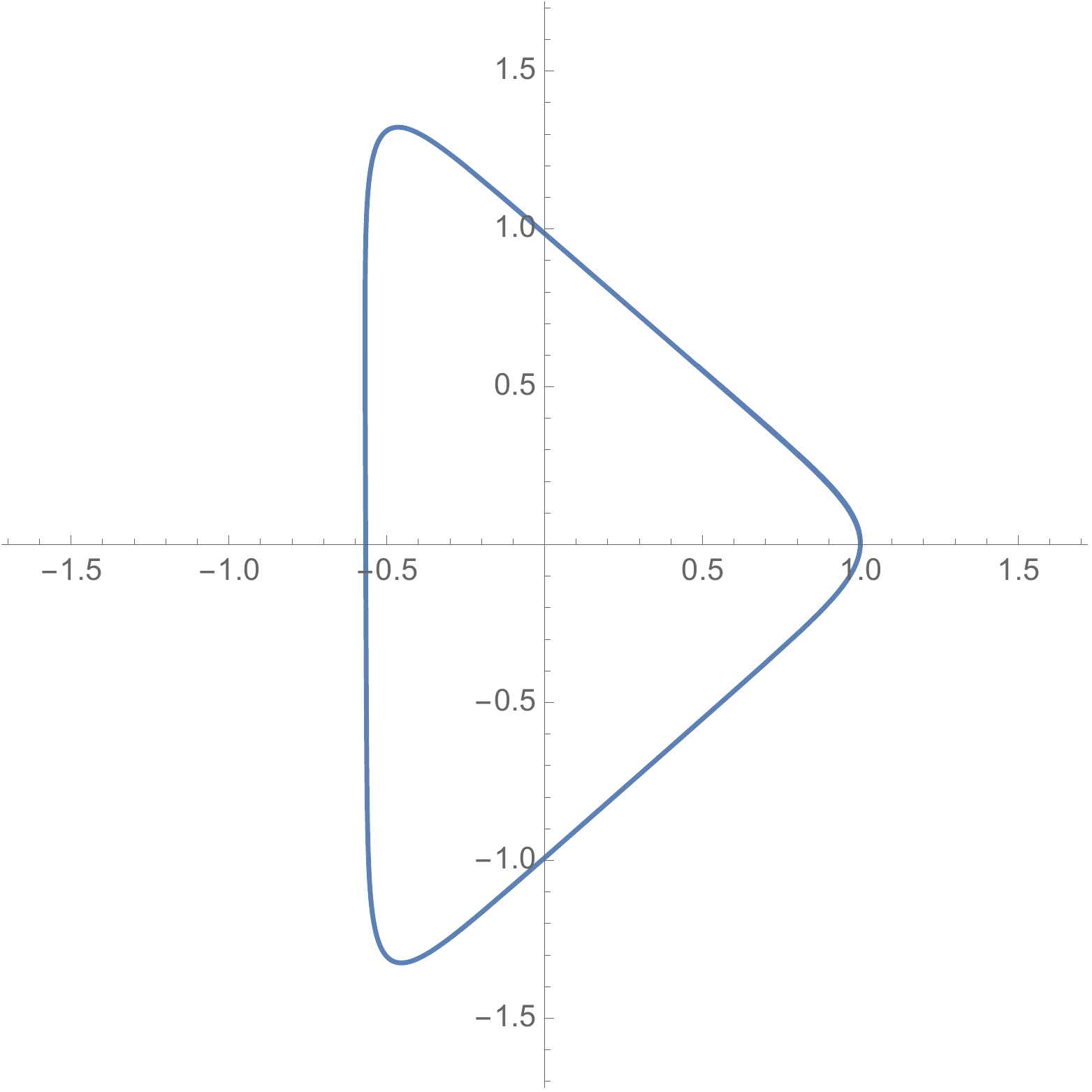}
\end{subfigure}
\begin{subfigure}{.512\textwidth}
\centering
\includegraphics[width=2.5in]{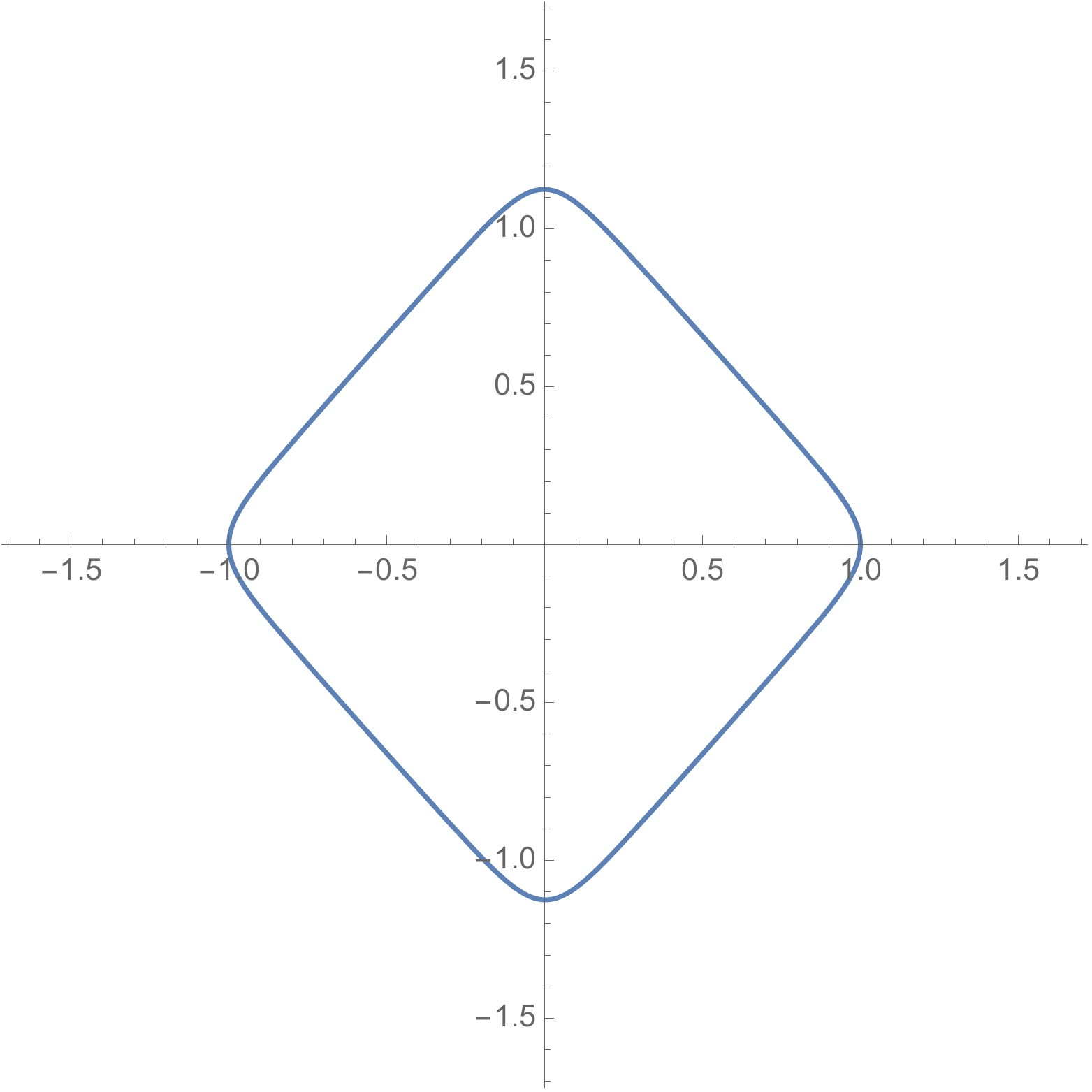}
\end{subfigure}
\begin{subfigure}{.512\textwidth}
\centering
\includegraphics[width=2.5in]{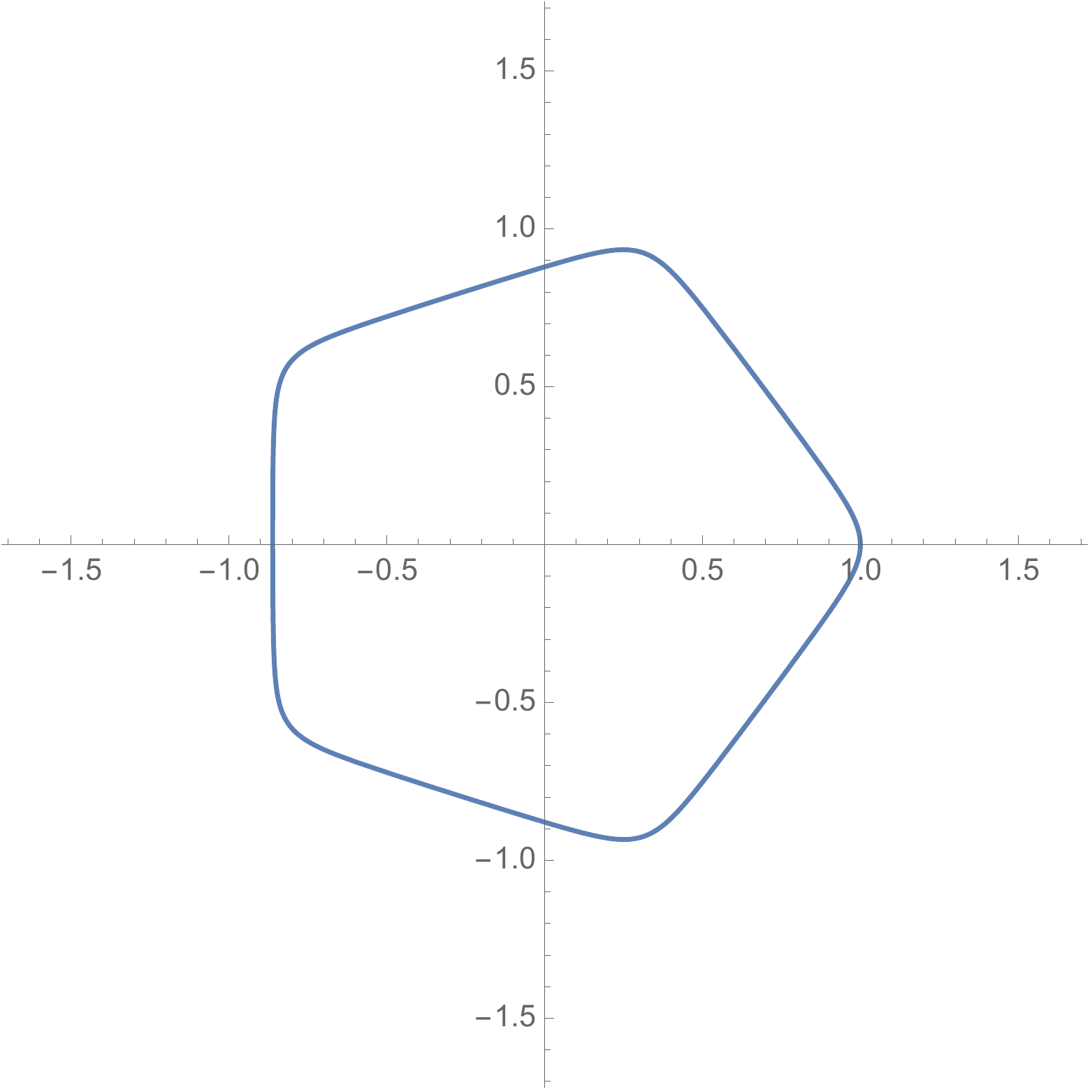}
\end{subfigure}
\begin{subfigure}{.5\textwidth}
\centering
\includegraphics[width=2.5in]{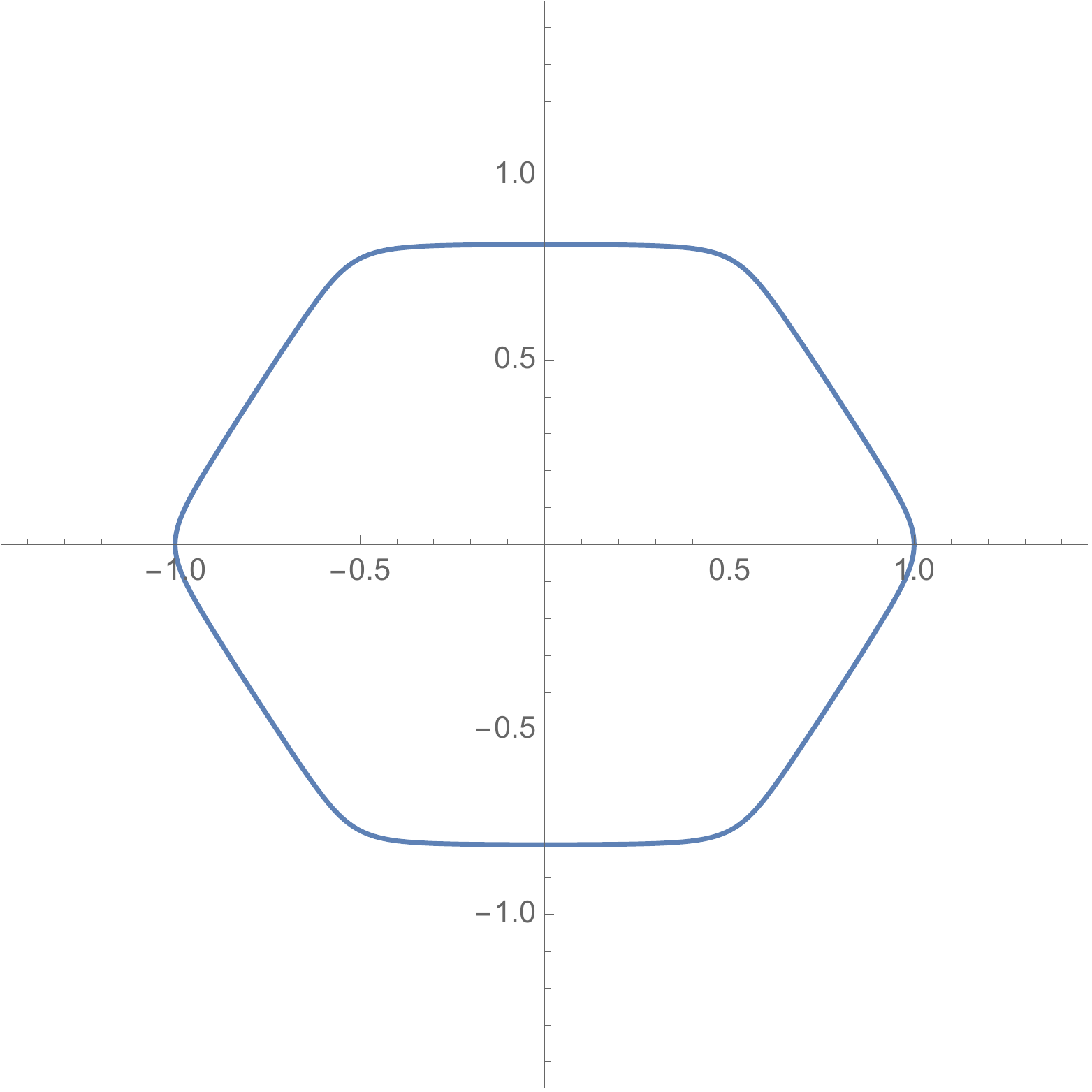}
\end{subfigure}
\begin{subfigure}{.512\textwidth}
\centering
\includegraphics[width=2.5in]{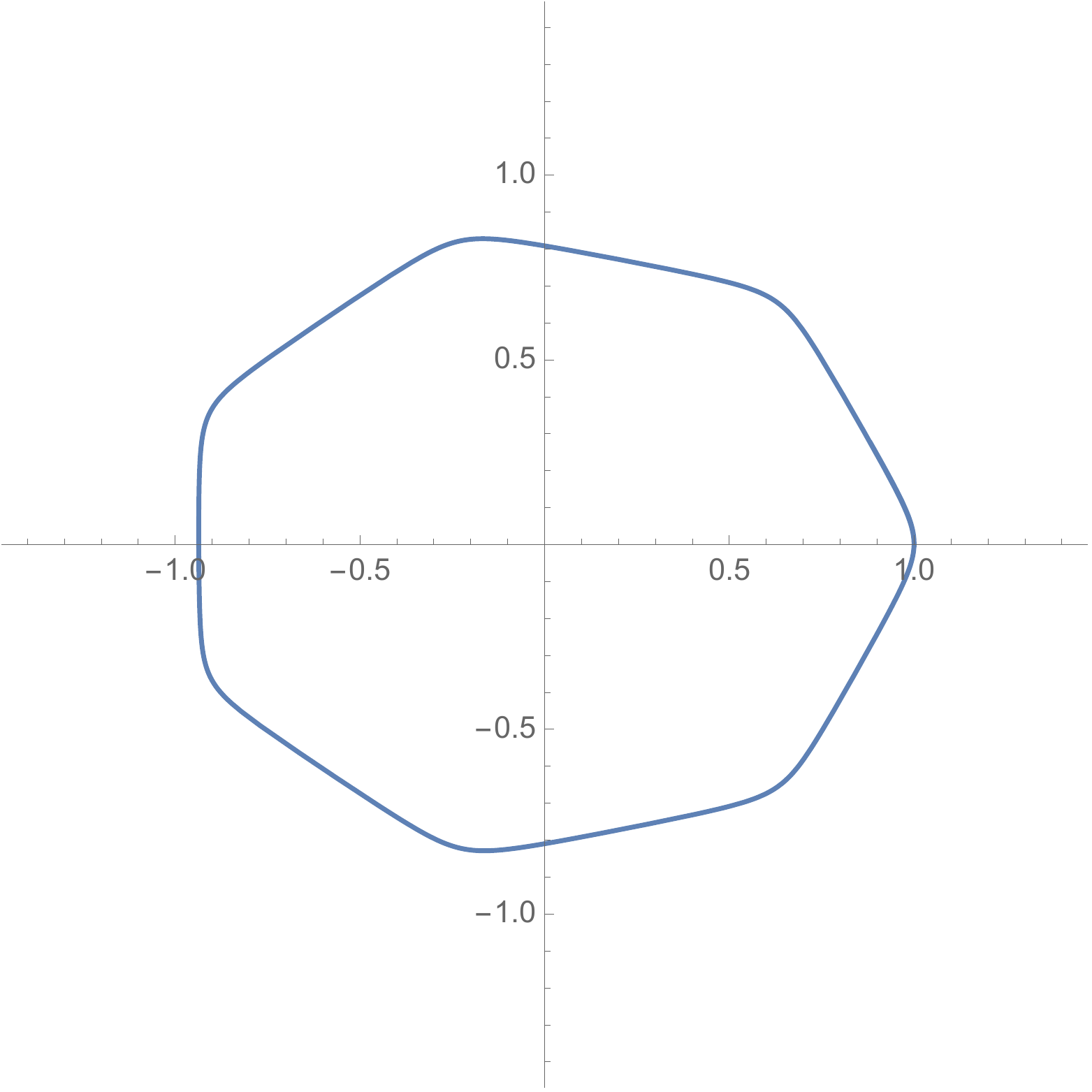}
\end{subfigure}
\begin{subfigure}{.5\textwidth}
\centering
\includegraphics[width=2.5in]{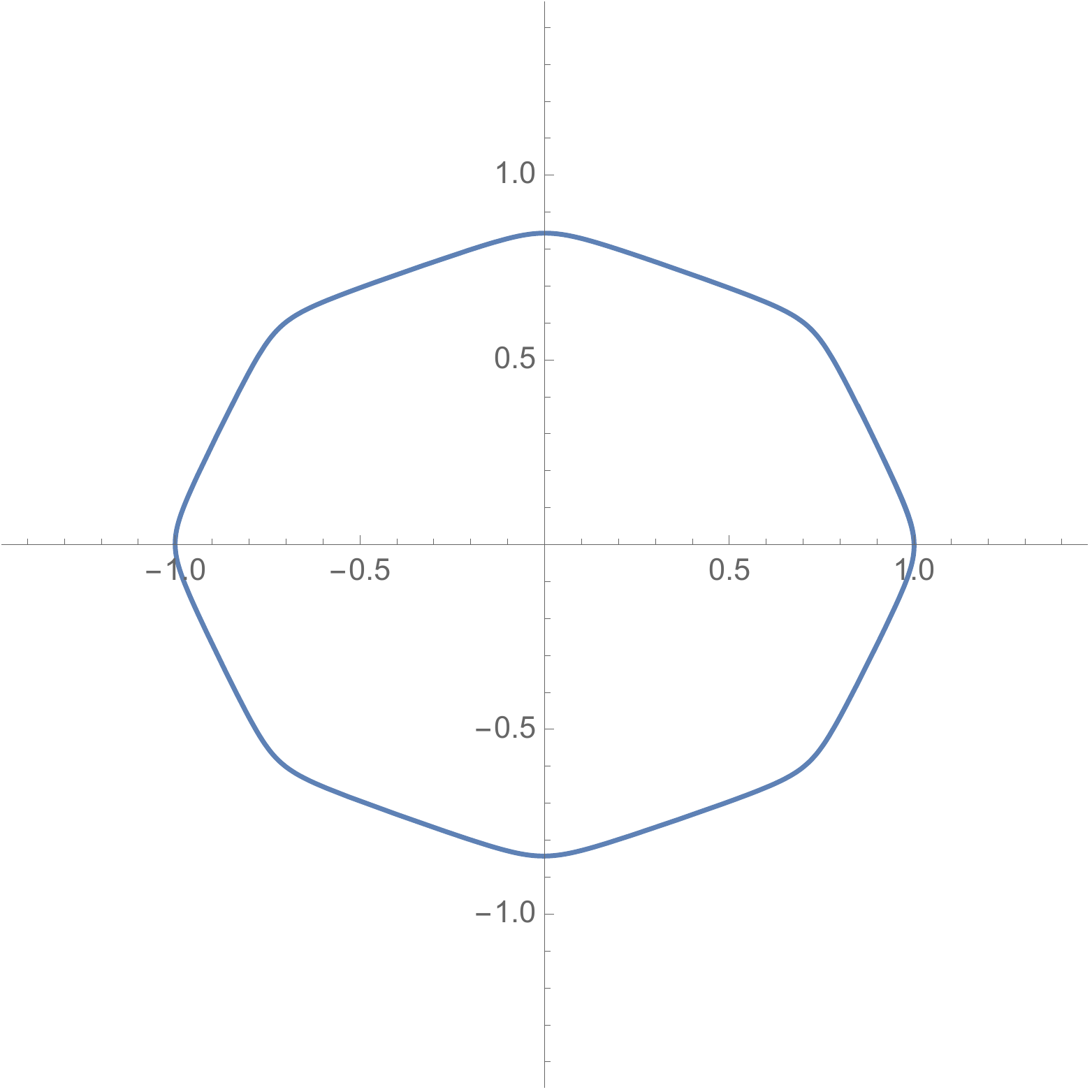}
\end{subfigure}
\caption{$a=\frac{7}{6}$, $a=\frac{23}{21}$, $a=\frac{17}{16}$, $a=\frac{47}{45}$, $a=\frac{31}{30}$, $a=\frac{79}{77}$}\label{fig:infinitesimal_deformations}
\end{figure}

A more systematic computer experiment (again using the Mathematica applet by Gil Bor) leads to Figure \ref{fig:infinitesimal_deformations} where curves corresponding to $k=4,5,6,7,8,9$ are displayed. We recall that for any value $a$, the functional $\B$ is invariant under $SL(2,\R)$. These correspond, of course, to trivial deformations. 

Figure \ref{fig:infinitesimal_deformations} leads us to the conjecture that for $a=\frac{k^2-2}{k^2-4}$ there exists an extremal curve which is a ``rounded $(k-1)$-gon''. Unfortunately, our software is currently not powerful enough to verify this. We hope to return to this point in the future. It is worth pointing out that these $(k-1)$-gons seem to approach a circle for $k$ very large. This is, at least, consistent with $a\to1$ when $k\to\infty$ and the circle is indeed rigid for $\mathcal{B}_1$. As a rather special case, the egg (case $a=1.4$ in Figure \ref{fig:egg}) should be considered a rounded 2-gon.

We collect some further questions.
\begin{itemize}  
\item What is the minimal rotation number of a periodic solution of $\B$ in dependence of $a$? It seems that this minimal rotation number goes to infinity with $a$.
\item In general, what happens if $a\to\infty$ respectively, $a\to-\infty$?
\item Is there some kind of duality for positive and negative values of $a$?
\item Is there an appropriate gradient flow of $\B$ similar to the metric case, see \cite{OPW} and \cite{Wan}?
\end{itemize}

\clearpage

\section{Appendix} \label{app}

In this appendix we describe a curious symmetry of equations (\ref{eq:1stord}). In the discovery of this symmetry, we were motivated by Bohlin's theorem as described in Appendix 1 of  \cite{ArB}.

Consider the family of equations
\begin{equation} \label{eq:3term}
\left(\frac{dF}{dt}\right)^2=uF^q+vF+w,
\end{equation}
where $F(t)$ is an unknown function and $u,v,w,q$ are parameters. We are looking for changes of independent and dependent variables
$$
F(t)=G(\tau)^\mu,\ \frac{d\tau}{dt}=G^\lambda
$$
that preserve the form of equation (\ref{eq:3term}), but possibly change the parameters $u,v,w$, and $q$. 

\begin{theorem} \label{thm:permut}
These changes of variables form a group, the group of permutations $S_3$. The orbit of the exponent $q$ is
$$
\left\{q, \frac{1}{q}, 1-q, \frac{1}{1-q}, \frac{q}{q-1},  \frac{q-1}{q}\right\}.
$$
\end{theorem}

We note that this is precisely how the permutations of four points in the projective line affect their cross-ratio.

\proof
Denote the new parameters by $\bar u,\bar v,\bar w$,  and $\bar q$.

Using the chain rule, one obtains the differential equation on $G$:
$$
\left(\frac{dG}{d\tau}\right)^2=
\frac{u}{\mu^2}G^{\mu q -2\mu-2\lambda+2}+\frac{v}{\mu^2}G^{2-2\lambda-\mu}+\frac{w}{\mu^2}G^{2-2\lambda-2\mu}.
$$
For this equation to have the same form as (\ref{eq:3term}), one needs the following relation between the exponents to hold:
$$
\{0,1\}\subset \{\mu q -2\mu-2\lambda+2,2-2\lambda-\mu,2-2\lambda-2\mu\}.
$$
Thus one needs to consider six cases. We present one of them: 
$$
\left\{
\begin{aligned}
&\mu q -2\mu-2\lambda+2=0\\
&2-2\lambda-\mu=1.
\end{aligned}
\right.
$$
Hence 
$$
\mu=\frac{1}{1-q},\ \bar q=2-2\lambda-2\mu =\frac{q}{q-1},\ \bar u=w(q-1)^2,\ \bar v= v(q-1)^2,\ \bar w=u(q-1)^2.
$$
The other five cases are analyzed in a similar way.
\proofend

Returning to equation (\ref{eq:1stord}), one has $q=b+2=\frac{2a-1}{a-1}$.
The $S_3$-orbit of the exponent $a$  is as follows:
$$
\left\{a, 1-a,  \frac{2a-1}{3a-1},  \frac{a}{3a-1}, \frac{2a-1}{3a-2}, \frac{a-1}{3a-2}\right\}.
$$
In particular, the orbit of 1 is $\{1, 1/2, 0\}$ and the orbit of 1/3 is $\{1/3, 2/3, \infty\}$, these are
special values in our study. All other  orbits consist of six elements.

We end with the following remarks. First, the function $G$ (or equivalently $F$) needs to be positive in order for $\frac{d\tau}{dt}=G^\lambda$ being an actual change of coordinates. This is always satisfied in our situation. Moreover, $F$ is periodic if and only if $G$ is. However, in the transition from a solution $F$ of \eqref{eq:1stord}, a special case of \eqref{eq:3term}, to a curve $\gamma$ it seems hard to see if $\g$ is again periodic. Therefore, the groups of coordinate changes above might or might not map closed curves to closed curves.



\begin{thebibliography}{99}


\bibitem{Ar} V. Arnold. {\it The geometry of spherical curves and quaternion algebra.} Russian Math. Surveys {\bf 50} (1995),  1--68.

\bibitem{ArB} V. Arnold. {\it Huygens and Barrow, Newton and Hooke. Pioneers in mathematical analysis and catastrophe theory from evolvents to quasicrystals.}  Birkh\"auser Verlag, Basel, 1990.

\bibitem{Ben} D. Bennequin. {\it Entrelacements et \'equations de Pfaff}. Third Schnepfenried geometry conference, Vol. 1 (Schnepfenried, 1982), 87--161,  Ast\'erisque, 107-108, Soc. Math. France, Paris, 1983. 

\bibitem{BBT} M. Bialy, G. Bor, S. Tabachnikov. {\it  Self-B\"acklund curves in centroaffine geometry and Lamé's equation}. arXiv:2010.02719.

\bibitem{La} R. Corless, G. Gonnet, D. Hare, D. Jeffrey, D. Knuth. {\it On the Lambert $W$ function}. Adv. Comput. Math. {\bf 5} (1996), 329--359.  

\bibitem{Lut} E. Lutwak. {\it Selected affine isoperimetric inequalities}. Handbook of convex geometry, Vol. A, 151--176, North-Holland, Amsterdam, 1993.

\bibitem{OPW} S. Okabe, P. Pozzi, G. Wheeler. {\it A gradient flow for the $p$-elastic energy defined on planar curves}. Math. Ann. {\bf 378} (2020), 777--828.

\bibitem{OT} V. Ovsienko, S. Tabachnikov. {\it Sturm theory, Ghys theorem on zeroes of the Schwarzian derivative and flattening of Legendrian curves.} Selecta Math. (N.S.) {\bf 2} (1996),  297--307.

\bibitem{Pin} U. Pinkall. {\it Hamiltonian flows on the space of star-shaped curves.} Results Math. {\bf 27} (1995), 328--332.

\bibitem{PS} G. P\'olya, G.  Szeg\"o. {\it Problems and theorems in analysis. II. Theory of functions, zeros, polynomials, determinants, number theory, geometry.}  Springer-Verlag, Berlin, 1998.

\bibitem{ST} G. Sapiro, A. Tannenbaum.  {\it On affine plane curve evolution.} J. Funct. Anal. {\bf 119} (1994), 79--120.

\bibitem{Sin} D. Singer. {\it Lectures on elastic curves and rods.} Curvature and variational modeling in physics and biophysics, 3--32, AIP Conf. Proc., 1002, Amer. Inst. Phys., Melville, NY, 2008. 

\bibitem{Wan} K. Wanatave. {\it Planar $p$-elastic curves and related generalized complete elliptic integrals.} Kodai Math. J. {\bf 37} (2014), 453--474.

\bibitem{WW} E. T. Whittaker, G. N. Watson. {\it A course of modern analysis. An introduction to the general theory of infinite processes and of analytic functions: with an account of the principal transcendental functions.} Fourth edition. Reprinted Cambridge University Press, New York, 1962.

\end{thebibliography}
\end{document}